\documentclass{article}%
\usepackage{amsmath}%
\usepackage{amsfonts}%
\usepackage{amssymb}%
\usepackage{graphicx}

\begin{document}

\title{Inversion Formula}
\author{Henrik Stenlund\thanks{The author is grateful to Visilab Signal Technologies for supporting this work.}
\\Visilab Signal Technologies Oy, Finland}
\date{July 27, 2010}
\maketitle

\begin{abstract}
This work introduces a new inversion formula for analytical functions. It is simple, generally applicable and straightforward to use both in hand calculations and for symbolic machine processing. It is easier to apply than the traditional Lagrange-B\"{u}rmann formula since no taking limits is required. This formula is important for inverting functions in physical and mathematical problems. The remainder term has been evaluated both for real and complex variable cases to give an error estimate for a broken series. \footnote{Visilab Report \#2010-07. Revision 1 written in \LaTeX with some sentences corrected, 2-1-2011. Revision 2 with equation (\ref{eqn90}) typo fixed. Revision 3 Dec. 2017 with derivation of the remainder term}
\subsection{Keywords}
inversion of functions, Taylor series, Lagrange-B\"{u}rmann inversion
formula, reversion of series
\subsection{Mathematical Classification}
Mathematics Subject Classification 2010: 11A25, 40E99, 32H02
\end{abstract}

\section{Introduction}
\subsection{General}
\noindent The inversion of an analytic function $f(z)$ with $z, u\in C$ 
\begin{equation}
f(z)=u
\end{equation}
is defined as
\begin{equation}
z=g(u)
\end{equation}
There is no general simple method known to determine $g(u)$ unless the variable $z$ can be readily solved from $f(z)$. Lagrange \cite{Lagrange1770} was the first to find a useful series expansion. B\"{u}rmann \cite{Burmann1799} and \cite{Hindenburg1798} generalized it to the Lagrange-B\"{u}rmann formula.

Good \cite{Good1960} extended the Lagrange-B\"{u}rmann formula to multiple variables. His formula is known as the Lagrange-Good formula and Hofbauer \cite{Hofbauer1979} supplied the proof. A number of investigations has been published over the Lagrange-B\"{u}rmann formula for various applications, like Zhao \cite{Zhao2004} and Merlini et al. \cite{Merlini2006}. Sokal \cite{Sokal2009} recently introduced a new generalization of the Lagrange-B\"{u}rmann formula. We first express the Lagrange-B\"{u}rmann inversion formula which is the present standard method for calculating the inverse. The new inversion formula is derived next.

\subsection{The Lagrange-B\"{u}rmann Inversion Formula}
Lagrange \cite{Lagrange1770} and B\"{u}rmann \cite{Burmann1799} introduced an inversion formula for a function $f(z)$ of a complex variable $z$.
\begin{equation}
f(z) = u
\end{equation}
with $f$ being analytic at some point $z_0$ and the first derivative at $z_0$, is required to be nonzero.
\begin{equation}
[{\frac{df(z)}{dz}}]_{z_{0}}\neq{0}
\end{equation} 
$f(z)$ has a value $u_0$ at $z_0$. The inverse function is $g(u)$
\begin{equation}
z=g(u)=g(f(z))
\end{equation} 
The Lagrange inversion formula or the Lagrange-B\"{u}rmann formula is a Taylor series as follows.
\begin{equation}
z=z_{0}+ \sum^{\infty}_{n=1}\frac{(u-u_{0})^{n}}{n!}[\lim_{z\rightarrow{z_0}}[\frac{d^{n-1}}{dz^{n-1}}(\frac{z-z_0}{f(z)-u_0})^{n}]]  \label{eqn10}
\end{equation} 
Proof of this formula can be found in \cite{Lagrange1770} and \cite{Burmann1799}. Taking limits in terms in equation  (\ref{eqn10}) usually requires lengthy calculations and a repeated use of L'Hospital's rule to get rid of the singularity. All terms belonging to a certain coefficient need to be kept together to determine the limit properly. This may be a very laborious task in hand calculations.
\section{The Inversion Formula}
Using the annotation of the preceding chapter, let
\begin{equation}
u=f(z) \  z,u \in C 
\end{equation}
and $f(z)$ be analytic over the interior of a circle
\begin{equation}
r=\left|z-z_0\right| 
\end{equation}
Let the inverse function $g(u)$ be analytic over the interior of a circle $R_0$ at $u_0$
\begin{equation}
R_0=\left|u-u_0\right| \label{eqn15}
\end{equation}
We have a Taylor series
\begin{equation}
z=z_{0}+ \sum^{\infty}_{n=1}\frac{(u-u_{0})^{n}}{n!}[\frac{d^{n}}{du^{n}}g(u)]_{u_{0}} \label{eqn20}
\end{equation} 
This series converges over the circle $R_1$ (as in equation (\ref{eqn15})). The equation (\ref{eqn20}) is very difficult to be used any further as such. Higher derivatives of $g(u)$ are requested and to get them, one would need the $g(u)$. We can use derivatives of $f(z)$ instead of $g(u)$. In order to circumvent the generation of progressively complicated terms, we proceed as follows. Differentiate equation (\ref{eqn30}) below.
\begin{equation}
z=g(u) \   \label{eqn30}
\end{equation}
to obtain
\begin{equation}
\frac{d}{dz}z=1=(\frac{d}{du}g)(\frac{d}{dz}u)=(\frac{d}{du}g(u))(\frac{d}{dz}f(z)) \   \label{eqn40}
\end{equation}
and solve it as
\begin{equation}
\frac{d}{du}g(u)=\frac{1}{(\frac{d}{dz}f(z))} \   \label{eqn50}
\end{equation}
Differentiate (\ref{eqn50}) further and solve it for
\begin{equation}
\frac{d^2}{du^2}g(u)=\frac{1}{(\frac{d}{dz}f(z))}\cdot[{\frac{d}{dz}\frac{1}{(\frac{d}{dz}f(z))}}] \   \label{eqn60}
\end{equation}
In the same manner the n'th derivative would be solved as
\begin{equation}
\frac{d^n}{du^n}g(u)=\frac{1}{(\frac{d}{dz}f(z))}\cdot[{\frac{d}{dz}\frac{1}{(\frac{d}{dz}f(z))}}\cdot[{\frac{d}{dz}\frac{1}{(\frac{d}{dz}f(z))}}\cdot\cdot\cdot[{\frac{d}{dz}\frac{1}{(\frac{d}{dz}f(z))}}\cdot[{\frac{d}{dz}\frac{1}{(\frac{d}{dz}f(z))}}]]]] \   \label{eqn70}
\end{equation}
having $n-1$ derivatives acting on the right side in addition to the bracketed derivatives acting on $f(z)$ alone. We can rearrange the brackets yielding
\begin{equation}
\frac{d^n}{du^n}g(u)=[\frac{1}{\frac{d}{dz}f(z)}\cdot{\frac{d}{dz}]^{n-1}\frac{1}{(\frac{d}{dz}f(z))}} \   \label{eqn80}
\end{equation}
The multiplying factor is a differential operator acting on all terms to the right containing any dependence on $z$. Placing this result to equation (\ref{eqn20}) yields the simplified inversion formula
\begin{equation}
z=z_{0}+\sum^{\infty}_{n=1}\frac{(u-u_{0})^{n}}{n!}[[\frac{1}{\frac{df(z)}{dz}}\cdot{\frac{d}{dz}]^{n-1}\frac{1}{(\frac{df(z)}{dz})}}]_{z_{0}} \label{eqn90}
\end{equation}
The necessary, but not sufficient, condition for the new inversion formula to converge is that the first derivative of $f(z)$ must be nonzero at $z_0$. The radius of convergence $R_1$ must be evaluated for each resulting series. If a singularity would appear at $z_0$, a translation to a nearby point should be made.

\section{The Remainder Term}
The remainder term expresses the error caused by breaking the series at N'th term. Thus it is important when broken series are used for approximations and the resulting error must be evaluated. In the following it is derived.

Assuming $f(z)$ is analytic (referring to equation (\ref{eqn90})) and by using equation (\ref{eqn80}) it is straightforward to prove by induction the following expression
\begin{equation}
z=z_0+\sum^{N}_{n=1}\frac{(u-u_0)^{n}}{n!}[(\frac{1}{f'(z)}\frac{d}{dz})^{n-1}\frac{1}{f'(z)}]_{z_0}+   \nonumber
\end{equation}
\begin{equation}
\frac{1}{N!}\int^{u}_{u_0}dt\cdot(u-t)^{N}\left[(\frac{1}{f'(z)}\frac{d}{dz})^{N}\frac{1}{f'(z)}\right]_{z=g(t),t=f(z)} \label{eqn8110}
\end{equation}
Here $z=g(u)$ is the inverse function of $u=f(z)$. 
\subsection{Real Case}
For the real variable case, by using the mean value theorem and integrating, we get
\begin{equation}
R_N(u)=\frac{(u-u_0)^{N+1}}{(N+1)!}\left[(\frac{1}{\frac{df(z)}{dz}}\frac{d}{dz})^{N+1}\frac{1}{f'(z)}\right]_{u=f(\zeta),z=g(u)=\zeta} \label{eqn8020}
\end{equation}
Here $\zeta$ is an intermediate value between $u_0$ and $u$. This result is analogous to the Lagrange remainder for real functions. 
\subsection{Complex Case}
For the complex variable case we obtain by using the corresponding mean value theorem IV by Curtiss \cite{Curtiss1907} 
\begin{equation}
\int^{z}_{a}dw{f(w)}\phi{(w)}=f(a)\int^{a+\theta{(z-a)}}_{a}dw\phi{(w)} \label{eqn8050}
\end{equation}
After integration we arrive at
\begin{equation}
R_N(u)=\frac{(1-(1-\theta)^{N+1})[(u-u_0)]^{N+1}}{(N+1)!}\left[(\frac{1}{\frac{df}{dz}}\frac{d}{dz})^{N}\frac{1}{f'(z)}\right]_{z_0} \label{eqn8140}
\end{equation}
Here $\theta$ has a value according to
\begin{equation}
\left|{\theta-1}\right|<1   \label{eqn8170}
\end{equation}
\section{Conclusions}
The equation (\ref{eqn90}) represents a simple alternative to the Lagrange-B\"{u}rmann formula (equation (\ref{eqn10})). The Lagrange-B\"{u}rmann formula requires taking limits and repeated use of L'Hospital's rule to remove the singularity. The new formula requires only elementary differentiation and evaluation at $z_0$. 

Comparison of coefficients in each term between the two formulas is not possible since the expansions are based on polynomials of $u$. A special case appears when $u_0=0$ making the expansions powers of $u$. This leads to equalities but not directly. One has to approach the limit $(z\rightarrow0)$ in equation (\ref{eqn10}) finally reaching terms identical with equation (\ref{eqn90}). Working in the opposite way is not possible.

In spite of its simplicity, this inversion formula can be applied generally. It can be used for inversion of functions and polynomials and for reversion of series. It is valid also for real variables. It is useful for estimating the behavior of the inverse function at some point with a few beginning terms. The radius of convergence needs to be studied for each new series.

The remainder terms both for real variable and complex variable cases have been developed. When estimating the behavior of the inverse function at some point with a few beginning terms, one would need the remainder term for estimating the resulting error.

\end{document}